\documentclass[a4paper, 10pt]{article}
\usepackage{mathrsfs}
\usepackage{amsmath}

\usepackage{amsthm}
\usepackage{array}
\usepackage{cases}
\usepackage{fullpage,authblk}
\usepackage{epsf,epsfig,amsfonts,amsgen,amsmath,amssymb,amstext,amsbsy,amsopn,amsthm,lineno,subfigure}
\usepackage{graphicx}
\usepackage{epstopdf}

\makeatletter
\def\th@plain{%
  \upshape
  \itshape 
}
\makeatother

\makeatletter
\renewenvironment{proof}[1][\proofname]{\par
  \pushQED{\qed}%
  \normalfont \topsep6\p@\@plus6\p@\relax
  \trivlist
  \item[\hskip\labelsep
        \bfseries
    #1\@addpunct{.}]\ignorespaces
}{%
  \popQED\endtrivlist\@endpefalse
}
\makeatother

\newtheorem{theorem}{Theorem}
\numberwithin{theorem}{section}
\newtheorem{lemma}{Lemma}

\numberwithin{corollary}{section}
\newtheorem{proposition}{Proposition}
\newtheorem{conjecture}{Conjecture}
\numberwithin{conjecture}{section}
\newtheorem*{conjecture*}{Conjecture}

\newtheorem{claim}{Claim}

\theoremstyle{definition}

\newtheorem{subclaim}{Subclaim}

\theoremstyle{remark}
\usepackage[varg]{txfonts}
\usepackage{graphicx, epsfig, subfigure}
\usepackage{enumerate} 
\usepackage[square, numbers, sort&compress]{natbib} 

\usepackage{hyperref}

\numberwithin{equation}{section}
\newtheorem{remark}{Remark}
\begin{document}
\title{Decomposing edge-colored graphs under color degree constraints}
\author[1,2]{ Ruonan Li \thanks{Supported by CSC (No.~201506290097); E-mail: liruonan@mail.nwpu.edu.cn}}
\author[3]{Shinya Fujita \thanks{Supported by JSPS KAKENHI (No.~15K04979); E-mail: fujita@yokohama-cu.ac.jp}}
\author[4]{Guanghui Wang\thanks{Supported by NSFC (No.~11471193, 11631014); E-mail:
ghwang@sdu.edu.cn}}
\affil[1]{Department of Applied Mathematics, Northwestern Polytechnical University, Xi'an, 710072, P.R.~China}
\affil[2]{ Faculty of EEMCS, University of Twente, P.O. Box 217, 7500 AE Enschede, The Netherlands}
\affil[3]{International College of Arts and Sciences, Yokohama City University, 22-2, Seto, Kanazawa-ku, Yokohama, 236-0027 Japan}
\affil[4]{ School of Mathematics, Shandong University, Jinan, 250100, P.R. China}

\renewcommand\Authands{ and }

\maketitle

\begin{abstract}
For an edge-colored graph $G$, the minimum color degree of $G$ means the minimum number of colors on edges which are adjacent to each vertex of $G$. We prove that if $G$ is an edge-colored graph with minimum color degree at least $5$ then $V(G)$ can be partitioned into two parts such that each part induces a subgraph with minimum color degree at least $2$. We show this theorem by proving a much stronger form. Moreover, we point out an important relationship between our theorem and Bermond-Thomassen's conjecture in digraphs.

\end{abstract}

\bigskip


\section{\bf Introduction}
When we try to solve a problem in dense graphs, decomposing a graph into two dense parts sometimes plays an important role in the proof argument. This is because one can apply an induction hypothesis to
one of the parts so as to obtain a partial configuration, and then use the other part to obtain a desired configuration.
Motivated by this natural strategy, many work has been done along this line, and now we have a variety of results in this partition problem. To name a few, Stiebitz \cite{Stiebitz} showed a nice theorem, which states that
every graph with minimum degree at least $a+b+1$ can be decomposed into two parts $A$ and $B$ such that $A$ has minimum degree
at least $a$ and $B$ has minimum degree at least $b$. We see that the bound $a+b+1$ is best possible by considering the complete graph of order $a+b+1$.
By the same example, Thomassen \cite{carsten3,carsten4} conjectured that every $(a+b+1)$-connected graph can be decomposed into two parts $A$ and $B$
in such a way that $A$ is $a$-connected and $B$ is $b$-connected.
It was shown by Thomassen himself \cite{tho1} that if $b \leq 2$, then
the conjecture is true. However, rather surprisingly, even for the case $b=3$ this conjecture is widely open until now.
Likewise, there are some other partition problems to find the partition $V(G)=A\cup B$ so that both $A$ and $B$ have a certain property, respectively. The digraph version of this problem was proposed at the Prague Midsummer Combinatorial Workshop in 1995: For a digraph $D$, let $\delta^+(D)$ be the minimum out degree of $D$.
For integers $s$ and $t$, does there exists a smallest value $f(s,t)$ such that each digraph $D$ with $\delta^+(D)\ge f(s,t)$ admits a
vertex partition ($D_1, D_2$) satisfying $\delta^+(D_1)\ge s$ and $\delta^+(D_2)\ge t$? In \cite{Alon: 1996, Alon2: 2006} Alon posed the problem: Is there a constant $c$ such that $f(1,2)\leq c$?We only know that $f(1,1)=3$ holds by a result of Thomassen~\cite{Thomassen: 1983}. No much progress has been made for this problem. Recently Stiebitz~\cite{St2} propose this problem again when he deals with the coloring number of graphs.
As observed from the above known results, it seems that these partition problems are very difficult even if we restrict our consideration to a very specific case.

In this paper, we would like to consider a similar problem in edge-colored graphs. To state our results, we introduce some notation and definitions. Throughout this paper, all graphs are finite and simple.
Let $G$ be an edge-colored graph. For an edge $e\in E(G)$, we use $col_G(e)$ to denote the color of $e$. For a vertex $v \in V(G)$, let $d^c_G(v)$ be the color degree of $v$ in $G$, that is, the number of colors on edges which are adjacent to $v$.
The minimum color degree of $G$ is denoted by $\delta^c(G) (:=\min\{d^c_G(v):\ v\in V(G)\})$.
For a subgraph $H$ of $G$ with $E(H)\neq\emptyset$, let $col_G(H)$ be the set of colors appeared in $E(H)$.
Also, for a pair of vertex-disjoint subgraphs $M, N$ in $G$, let $col_G(M,N)$ be the set of colors on edges between $M$ and $N$ in $G$.
For a vertex $v$ of $G$, let $N^c_G(v)=col_G(v, N_G(v))$. By definition, note that $d^c_G(v)=|N^c_G(v)|$. When there is no ambiguity, we often write $col(e)$ for $col_G(e)$, $col(H)$ for $col_G(H)$, $col(M,N)$ for $col_G(M,N)$ and $d^c(v)$ for $d^c_G(v)$. A graph is called a \textit{properly colored} graph (briefly, \textit{PC} graph) if no two adjacent edges have the same color. Let $a$ and $b$ be integers with $a\ge b\ge 1$. A pair $(A,B)$ is called \textit{(a,b)-feasible} if $A$ and $B$ are disjoint, non-empty subsets of $V(G)$ such that $\delta^c(G[A])\geq a$ and $\delta^c(G[B])\geq b$; in particular, if $G$ contains an $(a,b)$-feasible pair $(A,B)$ with $V(G)=A\cup B$ then we say that \textit{$G$ has an $(a,b)$-feasible partition}.

Again, motivated by the same complete graph having mutually distinct colored edges (that is, the rainbow $K_{a+b+1}$), we propose the following conjecture.

\begin{conjecture}\label{ourconj}
Let $a,b$ be integers with $a\ge b\ge 2$, and $G$ be an edge-colored graph with $\delta^c(G)\ge a+b+1$.
Then $G$ has an $(a,b)$-feasible partition.
\end{conjecture}

The main purpose of our paper is to give the solution of this conjecture for the case  $a=b=2$.

\begin{theorem}\label{main}
Conjecture~\ref{ourconj} is true for $a=b=2$.
\end{theorem}

To consider our problem, utilizing the structure of minimal subgraphs $H$ with $\delta^c(H)\ge 2$ will be very important. An edge-colored graph $G$ is \textit{$2$-colored} if $\delta^c(G)\geq 2$. Specifically, we say a graph $G$ is \textit{minimally $2$-colored} if $\delta^c(G)\ge 2$ holds but any proper subgraph $H$ of $G$ has minimum color degree less than $2$ in $H$. By definition, note that, every PC cycle is a minimally $2$-colored graph.
An edge-colored graph obtained from two disjoint cycles by joining a path is \textit{a generalized bowtie} (more briefly, call it \textit{g-bowtie}). We allow the case where the path joining two cycles is empty. In that case, the g-bowtie becomes a graph obtained from two disjoint cycles by identifying one vertex in each cycle.
Note also that $K_1+2K_2$ (that is, a graph obtained from two disjoint triangles by identifying one vertex of each triangle) is a g-bowtie with minimum order.

We have the following characterization of minimally $2$-colored graphs, which will be used to prove our main result.

\begin{theorem} \label{th4.1}  If an edge-colored graph $G$ is minimally $2$-colored, then $G$ is either a PC cycle or a 2-colored g-bowtie without containing PC cycles.
\end{theorem}

In fact Theorem~\ref{main} will be given by proving a much stronger result. We generalize the concept of $(a,b)$-feasible partitions as follows. For $k\ge 2$ if $V(G)$ can be partitioned into $k$ parts $A_1, A_2,\ldots , A_k$ such that $\delta^c(G[A_i])\ge a_i$ holds for each $1\le i\le k$ then we say that \textit{$G$ has an $(a_1,a_2, \ldots ,a_k)$-feasible partition}. In this paper, we will mainly focus on the case where $(a_1,a_2,\ldots ,a_k)=(2,2,\ldots ,2)$.
For simplicity, let us call \textit{$2^k$-feasible partition} in this special case (thus, $(2,2)$-feasible partitions are equivalent to $2^2$-feasible partitions).
To state our result, we shall introduce the following theorem, which is on the existence of vertex-disjoint directed cycles in digraphs.

\begin{theorem}[Thomassen \cite{Thomassen: 1983}]\label{Thm:f(k)_dicycle}
	For each natural number $k$ there exists a (smallest) number $f(k)$ such that every digraph $D$ with $\delta^+(D)\ge f(k)$ contains $k$ vertex-disjoint directed cycles.
\end{theorem}


Bermond and Thomassen \cite{Bermond-Thomassen: 1981} conjectured that $f(k)=2k-1$ and Alon \cite{Alon: 1996} showed  that $f(k)\le 64k$.

As above, for $k\geq 1$ let $f(k)$ be a function such that every directed graph $D$ satisfying $\delta^+(D)\geq f(k)$ contains $k$ disjoint directed cycles. Define a function $g(k)$ as follows.
\[
g(k)=
\begin{cases}
2,&\text{$k=1$;}\cr
\max\{f(k)+1,g(k-1)+3\},&\text{$k\geq 2$.}\cr
\end{cases}
\]

Our main result is following.

\begin{theorem}\label{main2}
	Let $G$ be an edge-colored graph with $\delta^c(G)\geq g(k)$. Then $G$ has a $2^k$-feasible partition.
\end{theorem}

We then focus on the case $b=2$ in Conjecture~\ref{ourconj}.
We obtained the following partial result.

\begin{theorem} \label{th1.1}  Let $a$ be an integer with $a\ge 2$, and let $K_n$ be an edge-colored complete graph of order $n$ with $\delta^c(K_n)\ge a+3$. Then $K_n$ has an $(a,2)$-feasible partition.
\end{theorem}

Also, in \cite{FLZ}, it is shown that any edge-colored complete bipartite graph $K_{m,n}$ with $\delta^c(K_{m,n})\ge 3$ contains a PC $C_4$. This yields the following.

\begin{theorem} \label{thflz}  If an edge-colored complete bipartite graph $K_{m,n}$ satisfies $\delta^c(K_{m,n})\ge a+2$, then $K_{m,n}$ admits an $(a,2)$-feasible partition.
\end{theorem}

Regarding Conjecture \ref{ourconj} in the general case, by using the probabilistic method, we get the following result.

\begin{theorem}\label{prob}
Let $a,b$ be integers with $a\geq b\geq 1$. If $G$ is an edge-colored graph with $|V(G)|=n$ and $\delta^c(G)\geq 2 ln n +4(a-1)$, then $G$ has an $(a,b)$-feasible partition.
\end{theorem}
Although our results might look a bit modest, proving Conjecture~\ref{ourconj} even for the case $b=2$ seems quite hard. This is because we could give a big improvement on the Alon's bound $``64k"$ if it is true.

\begin{theorem}\label{main3}
If Conjecture~\ref{ourconj} is true for $b=2$, then $f(k)\le 3k-1$.
\end{theorem}

In view of Theorem~\ref{main3}, it tells us that solving Conjecture~\ref{ourconj} completely seems a very difficult problem.

This paper is organized as follows. In Sections \ref{pre}, \ref{sec2} and \ref{sec4}, we give the proofs of Theorems \ref{th4.1}, \ref{main2} and \ref{prob}, respectively. In Section~\ref{a2}, we prove Theorems~\ref{th1.1} and~\ref{main3}. In particular, Theorem~\ref{main3} is obtained by a much stronger result (see Proposition~\ref{final} in Section~\ref{a2}).

\section{Proof of Theorem~\ref{th4.1}}\label{pre}

In order to prove this theorem, we first introduce a structural theorem characterizing edge-colored graphs without containing PC cycles.

\begin{theorem}[Grossman and H\"{a}ggkist \cite{Grossman:1983}, Yeo \cite{Yeo: 1997}]
\label{Yeo}
Let $G$ be an edge-colored graph containing no PC cycles. Then there is a vertex $z\in V(G)$ such that no component of $G-z$ is joint to $z$ with edges of more than one color.
\end{theorem}

\textit{Proof of Theorem~\ref{th4.1}}.

Let $G$ be a minimally 2-colored graph. If $G$ contains a subgraph $H$ which is a PC cycle or a 2-colored g-bowtie without containing PC cycles, then $G=H$ (otherwise, by deleting vertices in $V(G)\setminus V(H)$ or edges in $E(G)\setminus E(H)$, we obtain a smaller 2-colored graph). Hence, it is sufficient to prove that if $G$ contains no PC cycle, then $G$ contains a 2-colored g-bowtie. Apply Theorem~\ref{Yeo} to $G$. Since $G$ is minimally $2$-colored, we may assume that $G$ is connected and there is a vertex $z\in V(G)$ such that $G-z$ consists of two components $H_1$ and $H_2$ with all the edges between $z$ and $H_i$ has color $i$ for $i=1,2$.

Let $zx_1x_2\cdots x_p$ and $zy_1y_2\cdots y_q$, respectively, be longest PC paths in $G\backslash H_2$ and $G\backslash H_1$ starting from $z$. Set $x_0=z$ and $y_0=z$. Since $d^c_{G\backslash H_2}(x)\geq 2$ and $d^c_{G\backslash H_1}(y)\geq 2$ for arbitrary vertices $x\in V(H_1)$ and $y\in V(H_2)$, we have $p,q\geq 2$ and there exist vertices $x_i$ and $y_j$ for some $i,j$ with $0\leq i\leq p-2$ and $0\leq j\leq q-2$ such that $col(x_px_i)\neq col(x_{p-1}x_p)$ and $col(y_qy_j)\neq col(y_{q-1}y_q)$. Since $G$ contains no PC cycle, we have $col(x_px_i)=col(x_ix_{i+1})$ and  $col(y_qy_j)=col(y_jy_{j+1})$. Together, the path $x_ix_{i-1}\cdots x_1zy_1y_2\cdots y_j$ and cycles $x_ix_{i+1}\cdots x_px_i$ and $y_jy_{j+1}\cdots y_qy_j$ form a 2-colored g-bowtie.

The proof is complete.\hfill $\Box$\\

\section{Proof of Theorem~\ref{main2}}\label{sec2}

First we prove the following proposition.

\begin{proposition}\label{p1}
Let $G$ be an edge-colored graph with $\delta^c(G)\geq a+b-1$. If $G$ contains an $(a,b)$-feasible pair, then there exists an $(a,b)$-feasible partition of $G$.
\end{proposition}
\begin{proof}
Let $(A,B)$ be an $(a,b)$-feasible pair such that $A\cup B$ is maximal. If $(A,B)$ is not an $(a,b)$-feasible partition, then $A\cup B=V(G)\backslash S$ with $S\neq \emptyset$. Since $(A,B)$ is maximal, $(A,B\cup S)$ is not a feasible pair. Hence there exists a vertex $x$ in $S$ such that $d^c_{G[B\cup S]}(x)\leq b-1$. Recall that $d_G^c(x) \geq a+b-1$. So $d_{G[A\cup x]}^c(x)\geq a$. Thus $(A\cup x,B)$ is a feasible pair, which is a contradiction with the maximality of $(A,B)$. This proves that $(A,B)$ is an $(a,b)$-feasible partition of $G$.
\end{proof}

It is easy to check that the following proposition is also true.

\begin{proposition}\label{p1}
	Let $G$ be an edge-colored graph with $\delta^c(G)\geq \sum_{i=1}^{k}(a_i-1)+1$. If $G$ contains $k$ disjoint subgraphs $H_1,H_2,\ldots,H_k$ such that $\delta^c(H_i)\geq a_i$ for $i=1,2,\ldots,k$, then $G$ admits an $(a_1,a_2,\ldots, a_k)$-feasible partition.
\end{proposition}

In what follows, we will keep the above propositions in mind and use these facts as a matter of course.

\textit{Proof of Theorem~\ref{main2}}.





	We prove the theorem by contradiction. Let $G$ be a counterexample such that $G$ is chosen according to the following order of preferences.

	(i) $k$ is minimum;
	(ii) $|G|$ is minimum;
	(iii) $|E(G)|$ is minimum;
	(iv) $|col(G)|$ is maximum.

	By the choice of $G$, we know that $\delta^c(G)=g(k)$, $k\geq 2$ and $G$ contains no rainbow triangles. Let $S_v=\{u: d^c_{G-v}(u)=d^c_{G}(u)-1\}$. Then the following two claims obviously hold:
	\begin{claim}\label{cl:S_noempty}
		$S_v\neq \emptyset$ for all $v\in V(G)$.
	\end{claim}
	\begin{claim}\label{cl:uv}
		For each edge $uv\in E(G)$, either $u\in S_v$ or $v\in S_u$.
	\end{claim}	
	Now we prove the following claims.
	\begin{claim}
	\label{cl:star}
	For each color $i\in col(G)$, the subgraph $G_i$ induced by edges colored by $i$ is a star.	
	\end{claim}
	\begin{proof}
		By the choice of $G$, we know that $G$ contains no monochromatic triangles or monochromatic $P_3$'s. Thus for every color $i\in col(G)$, each component of $G_i$ is a star. If $G_i$ contains more than one component, then color one of the components with a color not in $col(G)$. Thus, we get a counterexample with more colors than $G$, which contradicts to the choice of $G$.
	\end{proof}
	\begin{claim}\label{cl:reduction}
	For $u,v\in V(G)$, if $u\in S_v$ and $v\not\in S_u$, then $S_u\cap N_G(v)\neq \emptyset$.
	\end{claim}
	\begin{proof}
		Suppose to the contrary that there exist vertices $u,v\in V(G)$ satisfying $u\in S_v$, $v\not\in S_u$ and $S_u\cap N_G(v)=\emptyset$. Then $col(vu)$ appears only once at $u$ and more than once at $v$. By Claim \ref{cl:star}, the color $col(vu)$ can only appear at $
		\{v\}\cup S_v$, particularly, not at $S_u$. Now we construct a colored graph $G'$ by deleting the vertex $u$ and adding edges $\{vx: x\in S_u\}$ to $G$ with all of them colored by $col(vu)$ (since $S_u\cap N_G(v)=\emptyset$, this is possible without resulting multi-edges).
		For each vertex $x\in V(G')\backslash S_u$, we have $d^c_{G'}(x)=d^c_{G}(x)$. For each vertex $y\in S_u$, we have $N^c_{G'}(y)\subseteq (N^c_{G}(y)\backslash col(uy))\cup col(vu)$. Since the color $col(vu)$ does not appear at $S_u$, we have $d^c_{G'}(y)=|N^c_{G'}(y)|=|N^c_{G}(y)|=d^c_G(y)$. This implies that $\delta^c(G')\geq\delta^c(G)=g(k)$. Note that $|G'|=|G|-1$. By the assumption of $G$, we know that $G'$ must admit a $2^k$-feasible partition. By Theorem \ref{th4.1}, $G'$ contains $k$ disjoint subgraphs $H_1, H_2,\ldots, H_k$ such that $H_i$ is either a PC cycle or a minimally 2-colored g-bowtie without containing PC cycles for $i=1,2,\ldots, k$. If $\bigcup_{i=1}^{k}E(H_i)\subseteq E(G)$, then  we can find a $2^k$-partition of $G$ as desired, a contradiction.  If $\bigcup_{i=1}^{k}E(H_i)\not\subseteq E(G)$, then all the edges in $T=(\bigcup_{i=1}^{k}E(H_i))\setminus E(G)$ form a monochromatic star with the vertex $v$ as a center. Thus, without loss of generality, assume that $T\subseteq E(H_1)$.
		
		Since $H_1$ is either a PC cycle or a minimally 2-colored g-bowtie without containing PC cycles, for each vertex $a\in H_1$ and each color $j\in col(H_1)$, the color $j$ appears at most $2$ times at $a$ in $H_1$. Thus we have $1\leq |T|\leq 2$.


    \begin{figure}[h]
    	\label{fig:T1T2}
    	\centering
   	\subfigure[$|T|=1$]{
   		\includegraphics[width=0.32\textwidth]{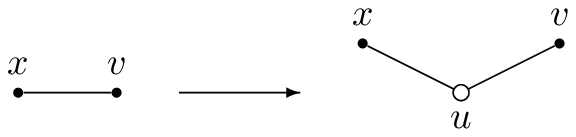}}
    	\hskip 70 pt
    	\subfigure[$|T|=2$]{
    	\includegraphics[width=0.32\textwidth]{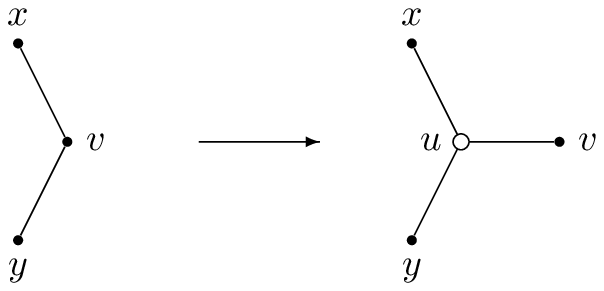}}
    	\caption{Cases of $|T|$}
    \end{figure}

		
		If $|T|=1$, then let $xv$ be the unique edge in $T$. Replace $xv$ in $H_1$ with the path $xuv$ (see Figure \ref{fig:T1T2}(a)). We obtain a colored graph $H'_1$ in $G$ with $\delta^c(H'_1)\geq 2$. Thus $H'_1,H_2,\ldots, H_k$ implies a $2^k$-feasible partition of $G$, a contradiction.
		
		If $|T|=2$, then let $T=\{vx,vy\}$. Since $col(vx)=col(vy)$, we know that $H_1$ is a minimally 2-colored g-bowtie with $v$ being an end vertex of the connecting path in $H_1$. Delete the edges $vx,vy$ and add vertex $u$ and edges $uv,ux,uy$ in $H_1$ (see Figure \ref{fig:T1T2}(b)). We obtain a g-bowtie $H'_1$ in $G$ with $\delta^c(H'_1)\geq 2$. Thus $H'_1,H_2,\ldots, H_k$ implies a $2^k$-feasible partition of $G$, a contradiction.
	\end{proof}
	\begin{claim}\label{cl:xy0}
		There exists an edge $xy\in E(G)$ such that $x\in S_y$ and $y\in S_x$.
	\end{claim}
	\begin{proof}
    Suppose not. Then by Claim \ref{cl:uv} , we can construct an oriented graph $D$ by orienting each edge $e=uv\in E(G)$ from $u$ to $v$ if and only if $v\in S_u$. Then $d^+_D(v)\geq 2$ for each vertex $v\in V(D)$. Let $T_i(v)=\{u: col(uv)=i\}$.
    \begin{subclaim}\label{scl:Ti}
    	For each vertex $v\in V(G)$ and colors $i,j\in col(G)$ with $i\neq j$, if $|T_i(v)|\geq 2$ and $|T_j(v)|\geq 2$, then the following statements hold:\\
    	$(a)$ $T_i(v)\cap T_j(v)=\emptyset$ and $E(T_i(v), T_j(v))=\emptyset$.\\ 	
    	$(b)$ $G[T_i(v)]$ contains at least one edge.
    \end{subclaim}
    \begin{proof}
    $(a)$	By the definition, we know that $T_i(v)\cap T_j(v)=\emptyset$. Since $|T_i(v)|\geq 2$ and $|T_j(v)|\geq 2$, we know that $T_i(v)\cup T_j(v)\subseteq S_v$. Let $u_i\in T_i(v)$ and $u_j\in T_j(v)$. Then colors $i$ and $j$ appears only once at $u_i$ and $u_j$, respectively. If $u_iu_j\in E(G)$, then $vu_iu_jv$ is a rainbow triangle, a contradiction. So we have $E(T_i(v), T_j(v))=\emptyset$.
    	    	
    $(b)$	Suppose that $G[T_i(v)]$ is empty for some color $i$ with $|T_i(v)|\geq 2$. Then choose $u\in T_i(v)$. We have $u\in S_v$ and $v\not\in S_u$. Apply Claim \ref{cl:reduction} to $u$ and $v$, we obtain $S_u\cap N_G(v)\neq\emptyset$. For each color $i'\in col(G)$ with $|T_{i'}(v)|\geq 2$, by Subclaim \ref{scl:Ti}$(a)$ and the assumption that $G[T_i(v)]$ is empty, we have $E(u,T_{i'}(v))=\emptyset$.
    Note that $$N^+_D(v)=\bigcup_{|T_{i'}(v)|\geq 2, i'\in col(G)} T_{i'}(v).$$
    We have $N_G(u)\cap N^+_D(v)=\emptyset$.
    Recall that $S_u\cap N_G(v)\neq \emptyset$ and $S_u\subseteq N_G(u)$. There must exist a vertex $x\in S_u \cap N_D^-(v) $. It is easy to check that $C=xuvx$ is a rainbow triangle in $G$, a contradiction.
    \end{proof}
    \begin{subclaim}\label{scl:one color}
    	For each vertex $v\in V(G)$, there is exactly one color $i\in col(G)$ such that $|T_i(v)|\geq 2$.
    \end{subclaim}
    \begin{proof}

        Given a vertex $v$, by Claim \ref{cl:S_noempty}, we can find a vertex $u\in S_v$. By the assumption of $G$, we have $v\not\in S_u$. Let $i=col(uv)$. Then $|T_i(v)|\geq 2$. This implies that for each vertex $v\in V(G)$, there is at least one color $i\in col(G)$ such that $|T_i(v)|\geq 2$.   Now, suppose to the contrary that there exists a vertex $v\in V(G)$ and colors $i,j\in col(G)$ with $i\neq j$ satisfying $|T_i(v)|\geq 2$ and $|T_j(v)|\geq 2$. By Subclaim \ref{scl:Ti}, we can choose edges $u_iw_i$ from $G[T_i(v)]$ and $u_jw_j$ from $G[T_j(v)]$. Let $F=G[v, u_i,w_i,u_j,w_j]$. Then $\delta^c(F)\geq 2$. Now we will discuss on the minimum color degree of $G-F$.
    	
    	If $\delta^c(G-F)\geq g(k-1)$, then by the assumption of $G$, $G-F$ has a $2^{k-1}$-feasible partition. Together with $G[V(F)]$, we obtain a $2^k$-feasible partition of $G$, a contradiction. So we have $\delta^c(G-F)< g(k-1)$. Let $x\in V(G-F)$ be a vertex satisfying $d^c_{G-F}(x)=\delta^c(G-F)$. Since $\delta^c(G)\geq g(k)\geq  g(k-1)+3$ and $|F|=5$, we have
    	$$4\leq |col(x,F)|\leq 5.$$
    	
    	For vertices $a\in \{u_i,w_i\}$ and $b\in  \{u_j,w_j\}$, if $|col(x,\{a,b,v\})|\geq 3$, then it is easy to check that either $xavx$  or $xbvx$ is a rainbow triangle, a contradiction. So we have $|col(x,\{a,b,v\})|\leq 2$. Note that $|col(x,F)|\geq 4$. This forces that $vx\not\in E(G)$ and $|col(x,\{u_i,w_i,u_j,w_j\})|=4$.
    	Thus $C=xu_ivu_jx$ is a rainbow cycle of length $4$. Suppose that there exists a vertex $y\in V(G-C)$ such that $d^c_{G-C}(y)< g(k-1)$. Then $|col(y,C)|\geq 4$. Note that $u_i,u_j\in S_v$. Thus either $yu_ivy$ or $yu_jvy$ is a rainbow triangle, a contradiction. Hence we have $\delta^c(G-C)\geq g(k-1)$. By the assumption of $G$, the graph $G-C$ has a $2^{k-1}$-feasible partition. Together with $G[V(C)]$, we get  a $2^{k}$-feasible partition of $G$, a contradiction.
    \end{proof}
    Subclaim \ref{scl:one color} implies that there are at least $g(k)-1$ colors appear only once at $v$ for each vertex $v\in V(G)$. Thus, we have $\delta^-(D)\geq g(k)-1\geq f(k)$. So $D$ contains $k$ disjoint directed cycles, which correspond to $k$ disjoint PC cycles in $G$, a contradiction.	
	\end{proof}
	
	\begin{claim}\label{cl:xy1}
		For each edge $xy\in E(G)$ satisfying $x\in S_y$ and $y\in S_x$, we have\\
		$(a)$ $|N^c_G(x)\cup N^c_G(y)-col(xy)|\leq g(k)-1$, and\\	
		$(b)$ $N_G(x)-y=N_G(y)-x=\{v_i:1\leq i\leq g(k)-1\}$, where $col(xv_i)=col(yv_i)$ and $col(xv_i)\neq col(xv_j)$ for $i,j\in[1,g(k)-1]$ with $i\neq j$.	
	\end{claim}
	\begin{proof}
	$(a)$ Since $G$ contains no rainbow triangles and $col(xy)$ appears only once at $x$ and $y$, respectively. we have $col(xu)=col(yu)$ for all $u\in N_G(x)\cap N_G(y)$. Now let $G'=G/xy$. Then $G'$ is well defined and $d^c_{G'}(v)=d^c_G(v)$ for all vertices in $V(G)\backslash\{x,y\}$. Let $z$ be the new vertex resulted by contracting the edge $xy$.
	
	Suppose that $|N^c_G(x)\cup N^c_G(y)-col(xy)|\geq g(k)$, then $d^c_{G'}(z)\geq g(k)$. Thus we have $\delta^c(G')\geq g(k)$. By the choice of $G$, we know that $G'$ must admit a $2^k$-feasible partition. By Theorem \ref{th4.1}, $G'$ contains $k$ disjoint subgraphs $H_1, H_2,\ldots, H_k$ such that $H_i$ ($i=1,2,\ldots, k$) is either a PC cycle or a minimally 2-colored g-bowtie without containing PC cycles.
	
	If $z\not\in \bigcup_{i=1}^{k}V(H_i)$, then $H_1, H_2,\ldots, H_k$ are $k$-disjoint subgraphs of $G$. This implies a $2^k$-feasible partition of $G$, a contradiction. So we can assume that $z\in V(H_1)$. Apparently, $2\leq d_{H_1}(z)\leq 4$.

	If $d_{H_1}(z)=2$, then let $N_{H_1}(z)=\{u,v\}$ (see Figure \ref{fig:dz2}). If $u,v\in N_G(x)$, then replace $z$ with $x$. If $u\in  N_G(x)$ and $v\not\in N_G(x)$, then replace the path $uzv$ with $uxyv$. In all cases, we can transform $H_1$ into a graph $H'_1\subseteq G$ such that $\delta^c(H'_1)\geq 2$ and $V(H'_1)\cap V(H_i)=\emptyset$ for $i=2,3 \ldots,k$.	Thus $H'_1,H_2,\ldots, H_k$ imply the existence of a $2^k$-feasible partition of $G$, a contradiction.


    \begin{figure}[h]
    	\centering
    	\includegraphics[width=0.50\linewidth]{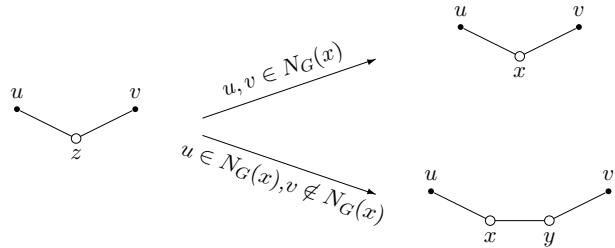}
    	\caption{$d_{H_1}(z)=2$}
    	\label{fig:dz2}
    \end{figure}
    \begin{figure}[h]
    	\centering
    	\includegraphics[width=0.50\linewidth]{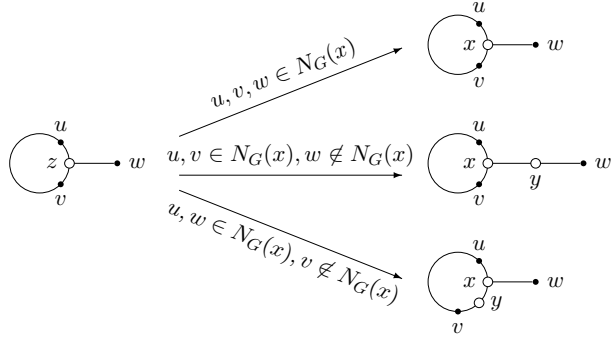}
    	\caption{$d_{H_1}(z)=3$}
    	\label{fig:dz3}
    \end{figure}
    \begin{figure}[h]
    	\centering
    	\includegraphics[width=0.65\linewidth]{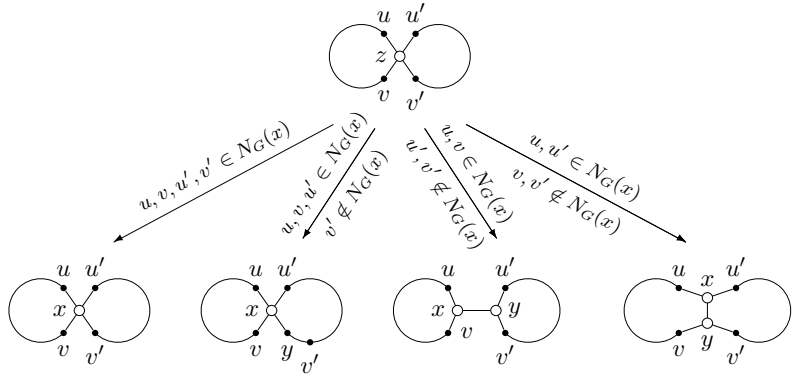}
    	\caption{$d_{H_1}(z)=4$}
    	\label{fig:dz4}
    \end{figure}

	If $d_{H_1}(z)=3$, then $H_1$ must be a minimally 2-colored g-bowtie with $z$ being an end-vertex of the connecting path. Let $N_{H_1}(z)=\{u,v,w\}$ with $u,v$ on a same cycle in $H_1$ (see Figure \ref{fig:dz3}). If $\{u,v,w\}\subseteq N_G(x)$, then replace $z$ with $x$. If $\{u,v\}\subseteq N_G(x)$ and $w\not\in N_G(x)$, then replace $zw$ with $xyw$. If $\{u,w\}\subseteq N_G(x)$ and $v\not\in N_G(x)$, then replace $zv$ with $xyv$. Constructions of the remaining cases are similar. Finally, in all cases, we can transform $H_1$ into a graph $H'_1\subseteq G$ such that $\delta^c(H'_1)\geq 2$ and $V(H'_1)\cap V(H_i)=\emptyset$ for $i=2,3 \ldots,k$.	Thus $H'_1,H_2,\ldots, H_k$ implies a $2^k$-feasible partition of $G$, a contradiction.
	
	If $d_{H_1}(z)=4$, then $H_1$ is a minimally 2-colored g-bowtie with two cycles overlapped on the vertex $z$. Let $N_{H_1}(z)=\{u,v,u',v'\}$ with $u,v$ on one cycle and $u',v'$ on the other cycle (see Figure \ref{fig:dz4}). If $\{u,v,u',v'\}\subseteq N_G(x)$, then replace $z$ with $x$.  If $\{u,v,u'\}\subseteq N^c(x)$ and $v'\not\in N^c(x)$, then replace the path $zv'$ with $xyv'$. If $\{u,v\}\subseteq N^c(x)$ and $\{u',v'\}\cap N^c(x)=\emptyset$, then split $z$ into the edge $xy$ such that the resulting graph is still a g-bowtie. If $\{u,u'\}\subseteq N^c(x)$ and $\{v,v'\}\cap N^c(x)=\emptyset$, then split $z$ into the edge $xy$ in an orthogonal direction such that the resulting graph is a cycle with one chord $xy$. Constructions of the remaining cases are similar. Finally, in all cases, we can transform $H_1$ into a graph $H'_1\subseteq G$ such that $\delta^c(H'_1)\geq 2$ and $V(H'_1)\cap V(H_i)=\emptyset$ for $i=2,3 \ldots,k$.	Thus $H'_1,H_2,\ldots, H_k$ implies a $2^k$-feasible partition of $G$, a contradiction.

	$(b)$	By Claim \ref{cl:xy1}$(a)$ and the fact that $d^c_G(x),d^c_G(y)\geq g(k)$, we have $N^c_G(x)=N^c_G(y)$ and  $d^c_G(x)=d^c_G(y)=g(k)$. For each color $j\in N^c_G(x)$ and $j\neq col(xy)$, since $G_{j}$ is a star and the color $j$ appears at $x$ and $y$, we know that $x,y$ must be leaf vertices of $G_j$. Let $v_j$ be the center of $G_j$. The proof is complete.
	\end{proof}
	
	Now let $\{x,y\}\cup \{v_i:1\leq i\leq g(k)-1\}$ be the set of vertices described in Claim \ref{cl:xy1}. Without loss of generality, let $col(xv_i)=i$ for $i\in[1,g(k)-1]$. Let $H$ be the subgraph of $G$ induced by $\{x,y\}\cup \{v_i:1\leq i\leq g(k)-1\}$ and $R=G-H$.
	\begin{claim}\label{cl:vi}
	For $1\leq i\leq g(k)-1$, $col(v_i,S_{v_i})=\{i\}$.	
	\end{claim}
	\begin{proof}
		Suppose to the contrary that there exists a vertex $u\in S_{v_i}$ such that $col(uv_i)\neq i$.	
		If $u=v_j$ for some $j$ with $1\leq j\leq g(k)-1$ and $j\neq i$, then $col(uv_i)=j$ (since $xv_iv_jx$ is not a rainbow triangle). Since the color $j$ appears at least 2 times at $v_j(=u)$, we know that $u\not\in S_{v_i}$, a contradiction. Now the vertex $u$ must belong to $V(R)$. Since each $G_j$ $(1\leq j\leq g(k)-1)$ is a star and $col(uv_i)\neq i$, we have $col(uv_i)\not\in [1,g(k)-1]$.
        If $v_i\in S_u$, then by applying Claim \ref{cl:xy1} to the edge $uv_i$, we have $N_G(u)-v_i=N_G(v_i)-u$. Since $x\in N_G(v_i)$, we have $x\in N_G(u)$, namely, $u\in N_G(x)$, a contradiction. So we have $v_i\not\in S_u$. Applying Claim \ref{cl:reduction} to $uv_i$, we obtain a vertex $v\in S_u\cap N_G(v_i)$. Note that $col(uv_i)\not\in [1,g(k)-1]$ and $G$ contains no rainbow triangle, we have $v\in R-u$. Let $F=G[x,y,v_i,u,v]$. It is easy to check that $\delta^c(F)\geq 2$.

        We will show that for each vertex $z\in G-F$, $|col(z,F)|\leq 3$. For $z\in R\cap (G-F)$, the assertion holds since $z$ has no neighbor to $x$ or $y$. Thus we may assume that $z=v_j$ for some $j$ with $1\leq j\leq g(k)-1$ and $j\neq i$.
        If $zv_i \notin E(G)$ or $col(zv_i)=j$, then we have the desired conclusion. So we may assume that
        $z$ is adjacent to $v_i$ and $col(zv_i)=i$ (otherwise, $zxv_iz$ is a rainbow triangle). Since there is no rainbow triangle and $G_i$ is a star, we can easily check that $zu\not\in E(G)$. So $z$ satisfies the desired property.

        Now, $\delta^c(G-F)\geq g(k)-3\geq g(k-1)$. So $G-F$ admits a $2^{k-1}$-feasible partition. Together with $G[V(F)]$, we obtain a $2^k$-feasible partition of $G$, a contradiction.
	\end{proof}
	
	\begin{claim}\label{cl:onlyxy}
	There exists a vertex $v_i$ with $1\leq i\leq g(k)-1$ such that $S_{v_i}=\{x,y\}$.	
	\end{claim}
	\begin{proof}
		Suppose not. Then there exists a vertex $u_i\in S_{v_i}\backslash \{x,y\}$ for all $i$ with $1\leq i\leq g(k)-1$. By Claim \ref{cl:vi}, $col(u_iv_i)=i$ for $1\leq i\leq g(k)-1$. Let $G'=G-\{x,y\}$. Then $\delta^c(G')\geq \delta^c(G)\geq g(k)$. By the choice of $G$, the graph $G'$ must admit a $2^k$-feasible partition, which implies that $G$ has a $2^k$-feasible partition, a contradiction.	
	\end{proof}
	
	We are now in a position to prove the theorem. Let $v_i$ be the vertex in Claim \ref{cl:onlyxy}. Since $d^c_H(v_i)\leq g(k)-1$ and $d^c_G(v_i)\geq g(k)$, there is a vertex $u\in R\cap N_G(v_i)$. Note that $u\not\in S_{v_i}$. By Claim \ref{cl:uv}, we have $v_i\in S_u$. Now apply Claim \ref{cl:reduction} to the edge $uv_i$, we have $S_{v_i} \cap N_G(u)\neq \emptyset$. This implies that either $x\in N_G(u)$ or $y\in N_G(u)$, a contradiction.
	
	This completes the proof of Theorem~\ref{main2}.\hfill $\Box$
\section{Proof of Theorem \ref{prob}}
\label{sec4}
\begin{lemma}\label{lem:ball}
	Let $k,x_1,x_2,\ldots,x_k$ be positive integers and $x_0$ a non-negative integer with $0\leq x_0\leq\frac{k}{2}$. Let $\{v_i^j: 1\leq i\leq k, 1\leq j\leq x_i\}$ be a set of $\sum_{i=1}^{k}x_i$ vertices such that each vertex $v_i^j$ is colored by $i$. Divide these vertices into two sets $S$ and $T$, randomly and independently, with $Pr(v_i^j\in S)=Pr(v_i^j\in T)=\frac{1}{2}$.  Let $P_S(x_0,x_1,\ldots,x_k)$ be the probability of the event that there are at most $x_0$ ($0\leq x_0\leq \frac{k}{2}$) differently colored vertices in $S$. Then
\begin{equation}\label{eq:prob}
	P_S(x_0,x_1,\ldots,x_k)\leq \sum_{j=0}^{x_0}\binom{k}{j}(\frac{1}{2})^k.
\end{equation}

\end{lemma}
\begin{proof}
	For convenience, we say a vector $\overrightarrow{x}=(x_0,x_1,x_2,\ldots,x_k)$ is {\it good} if $k, x_1,x_2,\ldots,x_k$ are positive integers and $x_0$ is a non-negative integer with $0\leq x_0\leq\frac{k}{2}$. Proving Lemma \ref{lem:ball} is equivalent to verify Inequation (\ref{eq:prob}) for all good vectors. For good vectors $\overrightarrow{x}=(x_0,x_1,\ldots,x_k)$ and $\overrightarrow{y}=(y_0,y_1,\cdots,y_{k'})$, we say $\overrightarrow{x}<\overrightarrow{y}$ if $(a)$ or $(b)$ holds.\\
	$(a)$ $k<k'$.\\
	$(b)$ $k=k'$ and there exists $t\in[1,k]$ such that  $x_t<y_t$ and $x_i=y_i$ for all $i$ with $0\leq i<t$.\\
    Now we will prove Inequation (\ref{eq:prob}) for every good vector $\overrightarrow{x}=(x_0,x_1,\ldots,x_k)$.
	
	By induction. First, it is easy to check that Inequation (\ref{eq:prob}) holds in the following three cases:
	$(1)$ $x_0=0$; 	$(2)$ $k=1$; $(3)$ $x_i=1$ for all $i$ with $1\leq i\leq k$. Now assume that $x_0\geq 1$, $k\geq 2$, $x_i\geq 2$ for some $i$ with $1\leq i\leq k$, and each good vector $\overrightarrow {y}$ with $\overrightarrow {y}<\overrightarrow {x}$ satisfies Inequation (\ref{eq:prob}). Consider the vertex $v_i^{1}$. We have
    $$
    P_S(\overrightarrow{x})=Pr(v_i^{1}\in T)P_S(x_0,x_1,\ldots,x_{i-1}, x_{i}-1,x_{i+1},\ldots,x_k)+Pr(v_i^{1}\in S)P_S(x_0-1,x_1,\ldots,x_{i-1},x_{i+1},\ldots,x_k).
    $$
    Let $\overrightarrow{y}=(x_0,x_1,\ldots,x_{i-1}, x_{i}-1,x_{i+1},\ldots,x_k)$ and $\overrightarrow{z}=(x_0-1,x_1,\ldots,x_{i-1},x_{i+1},\ldots,x_k)$. It is easy to see that $\overrightarrow{y}$ and $\overrightarrow{z}$ are good vectors with $\overrightarrow{y},\overrightarrow{z}<\overrightarrow{x}$. By induction hypothesis, we have
    $$P_S(\overrightarrow{y})\leq \sum_{j=0}^{x_0}\binom{k}{j}(\frac{1}{2})^k$$
    and
    $$P_S(\overrightarrow{z})\leq \sum_{j=0}^{x_0-1}\binom{k-1}{j}(\frac{1}{2})^{k-1}.$$
    Thus, we have
    \begin{eqnarray*}
    	P_S(\overrightarrow{x})
    	&\leq & \frac{1}{2}\sum_{j=0}^{x_0}\binom{k}{j}(\frac{1}{2})^k+\frac{1}{2}\sum_{j=0}^{x_0-1}\binom{k-1}{j}(\frac{1}{2})^{k-1}\\
    	&= &\frac{1}{2}\sum_{j=0}^{x_0}\binom{k}{j}(\frac{1}{2})^k+\sum_{j=1}^{x_0}\binom{k-1}{j-1}(\frac{1}{2})^{k}\\
     	&= &\frac{1}{2}\sum_{j=0}^{x_0}\binom{k}{j}(\frac{1}{2})^k+\sum_{j=1}^{x_0}\frac{j}{k}\binom{k}{j}(\frac{1}{2})^{k}\\
     	&\leq& \frac{1}{2}\sum_{j=0}^{x_0}\binom{k}{j}(\frac{1}{2})^k+\frac{x_0}{k}\sum_{j=1}^{x_0}\binom{k}{j}(\frac{1}{2})^{k}\\
     	&<& (\frac{1}{2}+\frac{x_0}{k})\sum_{j=0}^{x_0}\binom{k}{j}(\frac{1}{2})^{k}\\
     	&\leq&
        \sum_{j=0}^{x_0}\binom{k}{j}(\frac{1}{2})^{k}
    \end{eqnarray*}	
The proof is complete.
\end{proof}
{\it Proof of Theorem \ref{prob}}.
\begin{proof}
Assume $V(G)=\{1,2,\cdots,n\}$. We divide $V(G)$ into two disjoint parts $A,B$ randomly with $Pr(i\in A)=Pr(i\in B)=\frac{1}{2}$ for each vertex $i\in V(G)$. For each vertex $i\in A$,  the bad event $A_i$ means that for vertex $i$, $\{d_{G[A]}^c(i)\leq a-1\}$. By Lemma \ref{lem:ball}, we have
$$Pr(A_i)\leq\sum_{j= 0}^{a-1}\binom{d_G^c(i)}{j}(\frac{1}{2})^{d_G^c(i)}=\sum_{j= d_G^c(i)-a+1}^{d_G^c(i)}\binom{d_G^c(i)}{j}(\frac{1}{2})^{d_G^c(i)}
=Pr(X\geq d_G^c(i)-a+1),$$
where $X\sim B(d_G^c(i),\frac{1}{2})$.

Recall that Chernoff's bound: $Pr[X-E(X)\geq n\epsilon]<e^{-2n\epsilon^2}$, where $X\sim B(n,\frac{1}{2})$. We get
$$Pr(X\geq d_G^c(i)-a+1)=Pr(X-\frac{d_G^c(i)}{2}\geq \frac{d_G^c(i)}{2}-a+1)<e^{-2(\frac{d_G^c(i)}{2}-a+1)^2/d_G^c(i)}.$$
Since $d_G^c(i)\geq \delta^c(G)>2(a-1)$, we have
$$Pr(A_i)<e^{-2(\frac{d_G^c(i)}{2}-a+1)^2/d_G^c(i)}\leq e^{-2(\frac{\delta^c(G)}{2}-a+1)^2/{\delta^c(G)}}.$$ Similarly, for each vertex $j\in B$, the bad event $B_j$ means that $\{d_{G[B]}^c(j)\leq b-1\}$ and $Pr(B_j)<e^{-2(\frac{\delta^c(G)}{2}-b+1)^2/{\delta^c(G)}}\leq e^{-2(\frac{\delta^c(G)}{2}-a+1)^2/{\delta^c(G)}}$. So
$$Pr((\bigcup\limits_{i\in A} A_i)\cup(\bigcup \limits_{j\in B} B_j))\leq \sum\limits_{i\in A} Pr(A_i)+\sum\limits_{j\in B} Pr(B_j)< ne^{-2(\frac{\delta^c(G)}{2}-a+1)^2/{\delta^c(G)}}.$$
If $ne^{-2(\frac{\delta^c(G)}{2}-a+1)^2/{\delta^c(G)}}\leq 1$, which means $1-Pr[(\bigcup\limits_{i\in A} A_i)\cup(\bigcup \limits_{i\in B} B_i)]>0$, then $\frac{{\delta^c(G)}}{2}-2(a-1)+\frac{2(a-1)^2}{{\delta^c(G)}}\geq lnn$. The last inequality holds by the condition that ${\delta^c(G)}\geq 2lnn +4(a-1)$. Thus there exists a partition such that neither event $A_i$ nor $B_i$ happens. So we have an $(a,b)$-feasible partition.
\end{proof}

\section{From $(a,2)$-feasible partitions to Bermond-Thomassen's conjecture}\label{a2}

Firstly, we give the proof of Thorem~\ref{th1.1}.

\textit{Proof of Theorem~\ref{th1.1}.}

In order to prove the theorem, we use the following fact.

\begin{lemma}\label{L1}\cite{Gallai: 1976}
In any rainbow triangle-free coloring of a complete graph, there exists
a vertex partition $(V_1,V_2\ldots ,V_t)$ of $V(K_n)$ with $t\ge 2$ such that between the parts, there are a total of
at most two colors and, between every pair of parts $V_i, V_j$ with $i\neq j$, there is only one color on the edges.
\end{lemma}

If $K_n$ contains a rainbow triangle $C$, then let $A=C$ and $B=K_n-C$. It follows that $\delta^c(A)\geq 2$ and $\delta^c(B)\geq a$. So $(A,B)$ is an $(a,2)$-feasible partition. Now we assume that $K_n$ contains no rainbow triangle. Utilizing Lemma \ref{L1}, we can easily find an $(a,2)$-feasible partition. Thus Theorem~\ref{th1.1} holds.
\hfill $\Box$\\

In this section, we will point out a relationship between $(a,2)$-feasible partitions in edge-colored graphs and Bermond-Thomassen's conjecture in digraphs. In fact, Bermond-Thomassen's conjecture has not even been confirmed in multi-partite tournaments. Recently, Li et al. \cite{LR:2017+} revealed a relationship between PC cycles in edge-colored complete graphs and Bermond-Thomassen's conjecture on multi-partite tournaments.

We prove the following proposition.

	\begin{proposition}\label{prop:3relation}For $k\ge 1$ let $d_1,\ldots, d_k$ be positive integers, and
		let $f(d_1,d_2,\ldots,d_k)$, $g(d_1,d_2,\ldots,d_k)$ and $h(d_1,d_2,\ldots,d_k)$ be the minimum values which make the following three statements true:
\begin{description}
		\item{$(1)$} Every oriented graph $D$ with $\delta^+(D)\geq f(d_1,d_2,\ldots,d_k)$ has a vertex-partition $(V_1,V_2,
		\ldots,V_k)$ with $\delta^+(D[V_i])\geq d_i$ for $i=1,2,\ldots,k$.

		\item{$(2)$} Every edge-colored graph $G$ with $\delta^c(G)\geq g(d_1,d_2,\ldots,d_k)$ has a $(d_1,d_2,\ldots,d_k)$-feasible partition.
		\item{$(3)$} Every edge-colored complete graph $K$ with $\delta^c(K)\geq h(d_1,d_2,\ldots,d_k)$ has a $(d_1,d_2,\ldots,d_k)$-feasible partition.
\end{description}
		Then we have
		$$f(d_1-1,d_2-1,\ldots,d_k-1)\leq g(d_1,d_2,\ldots,d_k)\leq h(d_1+1,d_2+1,\ldots,d_k+1).$$
	\end{proposition}
	\begin{proof}
		Given an oriented graph $D$, we construct an edge-colored graph $G$ with $V(G)=V(D)$, $E(G)=\{uv: uv\in A(D) \text{~or~} vu\in A(D)\}$ and $col_G(uv)=v$ if and only if $uv\in A(D)$. If $\delta^+(D)\geq g(d_1,d_2,\cdots,d_k)$, then by the construction, we know that $\delta^c(G)\geq g(d_1,d_2,\cdots,d_k)$. Thus, $G$ admits a partition $V_1,V_2,\ldots,V_k$ such that $\delta^c(G[V_i])\geq d_i$ for $i=1,2,\ldots,k$. In turn, by the construction, we have $\delta^+(D[V_i])\geq d_i-1$ for $i=1,2,\ldots,k$. Recall the definition of function $f$. We know that
		$$f(d_1-1,d_2-1,\ldots,d_k-1)\leq g(d_1,d_2,\ldots,d_k).$$
		
		Given an edge-colored graph $G$, we construct an edge-colored complete graph $K$ with $V(K)=V(G)$, $col_K(e)=col_G(e)$ for all $e\in E(G)$, $col_K(e)=c_0$ for all $e\in E(K)\setminus E(G)$ and $c_0\not\in col(G)$. If $\delta^c(G)\geq h(d_1+1,d_2+1,\ldots,d_k+1)$, then $\delta^c(K)\geq h(d_1+1,d_2+1,\ldots,d_k+1)$. By the definition of $h$, we know that there exists a partition $V_1,V_2,\ldots,V_k$ of $K$ such that $\delta^c(K[V_i])\geq d_i+1$ for $i=1,2,\ldots,k$. By the construction of $K$, we have $\delta^c(G[V_i])\geq d_i$ for $i=1,2,\ldots,k$. Recall the definition of $g$. We know that
		$$g(d_1,d_2,\ldots,d_k)\leq h(d_1+1,d_2+1,\ldots,d_k+1).$$		
	\end{proof}
	\begin{remark}
		The existence of $f(d_1,d_2,\ldots,d_k)$ for $d_i\geq 2$ ($i=1,2,\ldots,k$) and $k\geq 2$ is still unknown according to  \cite{Alon: 1996}. Proposition \ref{prop:3relation} implies that we could show the existence of $f(d_1,d_2,\ldots,d_k)$ by proving the existence of $g(d_1+1,d_2+1,\ldots,d_k+1)$ or $h(d_1+2,d_2+2,\ldots,d_k+2)$.
	\end{remark}
When $d_1=d_2=\cdots=d_k=d$, for simplicity, we write $f(d,d,\cdots,d)_k$ instead of $f(d_1,d_2,\cdots,d_k)$. This also applies to functions $g$ and $h$.

The following result provides us the direct consequence of Theorem~\ref{main3}.

\begin{proposition}\label{final}
	If $g(a,2)\leq a+t$ for an integer $t$ and all $a\in \mathbb{N}$, then
	$$f(1,1,\ldots,1)_k\leq g(2,2,\ldots,2)_k\leq tk-t+2.$$
\end{proposition}
\begin{proof}
	According to Proposition \ref{prop:3relation}, we only need to prove that $g(2,2,\ldots,2)_k\leq tk-t+2$. By induction on $k$. Since $g(a,2)\leq a+t$ for all $a\in \mathbb{N}$. We have $g(2,2)\leq t+2$. Assume that $g(2,2,\ldots,2)_{k-1}\leq (k-2)t+2$. and let $x=g(2,2,\ldots,2)_{k-1}$. Then 	
	$$g(2,2,\ldots,2)_k\leq g(x,2)\leq x+t\leq (k-1)t+2=tk-t+2.$$
	So $g(2,2,\ldots,2)_k\leq tk-t+2$ for all $k\geq 2$.
	
	The proof is complete.
\end{proof}
\begin{remark}
    Bermond and Thomassen \cite{Bermond-Thomassen: 1981} conjectured that $f(1,1,\ldots,1)_k=2k-1$ (the conjecture is proposed for simple directed graphs and it is sufficient to prove it in oriented graphs). Recall that the best known upper bound of $f(1,1,\ldots,1)_k$ is $64k$ (by Alon \cite{Alon: 1996}). In view of Proposition~\ref{final}, we suggest that considering $(a,2)$-feasible partitions in edge-colored graphs could be a reasonable approach for improving Alon's result concerning Bermond-Thomassen's conjecture in digraphs.
\end{remark}

\end{document}